\documentclass[12pt]{article}
\usepackage{amsmath}
\usepackage{amsbsy}
\usepackage{amsfonts}
\usepackage{mathtools}
\usepackage{bbm}
\usepackage{amsmath}
\usepackage{amstext}
\usepackage{amsthm}
\usepackage{bbold}
\usepackage{amssymb}
\usepackage[dvips]{graphicx}
\usepackage{amsfonts}
\usepackage{caption2}
\usepackage{amsfonts,mathrsfs}
\usepackage{makeidx}
\usepackage{fancyhdr}
\usepackage{epic,eepic,epsfig}
\usepackage{mdwlist}
\usepackage{amscd}
\usepackage{amsbsy}
\usepackage{epsfig,graphicx}
\usepackage{epsfig,mathrsfs,wasysym,stmaryrd}
\usepackage[all]{xypic}
\usepackage[knot,poly]{xy}
\usepackage{latexsym}
\usepackage{comment}

\theoremstyle{plain}

\newtheorem{theorem}{Theorem}
\newtheorem{proposition}{Proposition}

\newtheorem{lemma}{Lemma}

\newenvironment{pf}{\medskip\noindent{Proof:}
  \hspace{-.5cm}      \enspace}{\hfill \qed \newline \smallskip}

\date{}

\begin{document}

\begin{center}
\textbf{\LARGE{ On rationality of generating function for the number of spanning trees in circulant graphs}}
\vspace{12pt}

{\large\textbf{A.~D.~Mednykh,}}\footnote{{\small\em Sobolev Institute of Mathematics,
Novosibirsk State University, smedn@mail.ru}}
{\large\textbf{I.~A.~Mednykh,}}\footnote{{\small\em Sobolev Institute of Mathematics,
Novosibirsk State University, ilyamednykh@mail.ru}}
\end{center}

\section*{Abstract}

Let $F(x)=\sum\limits_{n=1}^\infty\tau(n)x^n$ be the generating function for the number
$\tau(n)$ of spanning trees in the   circulant graphs $C_{n}(s_1,s_2,\ldots,s_k).$  We show that
$F(x)$ is a rational function with integer coefficients satisfying the property $F(x)=F(1/x).$ A similar result is also true for
the circulant graphs of odd valency  $C_{2n}(s_1,s_2,\ldots,s_k,n).$  We illustrate the obtained results by a series of examples.

\bigskip

\noindent
\textbf{Key Words:} spanning tree, circulant graph, Chebyshev polynomial, generating function\\
\textbf{AMS classification:} 05C30, 39A10\\

\section{Introduction}

The \textit{complexity} of a finite connected graph $G$, denoted by $\tau(G),$ is the number
of spanning trees of $G.$ One of the first results on the complexity was obtained by Cayley \cite{Cay89}
who proved that the number of spanning trees in the complete graph $K_n$ on $n$ vertices is $n^{n-2}.$

The famous Kirchhoff's Matrix Tree Theorem~\cite{Kir47} states that $\tau(G)$ can be expressed
as the product of nonzero Laplacian eigenvalues of $G$ divided by the number of its vertices.
Since then, a lot of papers devoted to the complexity of various classes of graphs were published.
In particular, explicit formulae were obtained for complete multipartite graphs~\cite{Cay89, Austin60},
almost complete graphs~\cite{Wein58}, wheels~\cite{BoePro}, fans~\cite{Hil74}, prisms~\cite{BB87},
ladders~\cite{Sed69}, M\"obius ladders~\cite{Sed70}, lattices \cite{Wu77, SW00, Louis}, anti-prisms \cite{SWZ16},
complete prisms~\cite{Sch74} and for many other families.

Starting with Boesch and Prodinger \cite{BoePro} the idea to study the complexity of graphs by making
use of Chebyshev polynomials was implemented. This idea provided a way to find complexity for different families of circulant
graphs and their natural generalisations  \cite{KwonMedMed, Med1, Xiebin, XiebinLinZhang, YTA97,
ZhangYongGol, ZhangYongGolin, Louis}.

The asymptotic behavior of complexity for some families of graphs can be investigated from
the point of view of so called Malher measure \cite{GutRog}, \cite{SilWil}, \cite{SilWil1}.
Mahler measure of a polynomial $P(z)$, with complex coefficients, is the product of the roots
of $P(z)$ whose modulus is greater than $1$ multiplied by the leading coefficient.

In the recent papers by  the authors \cite{MedMed3} and \cite{MedMed2018}, explicit  formulas for the number of spanning trees
$\tau(n)$ in circulant graphs $C_{n}(s_1,s_2,\ldots,s_k)$ and  $C_{2n}(s_1,s_2,\ldots,s_k,n)$  were obtained.
It was shown that in both cases the number of spanning trees can be represented in the form $\tau(n) = p\,n\,a(n)^2,$ where $a(n)$ is an integer sequence and $p$ is a prescribed natural number depending on the parity of $n.$ Also,   asymptotic formulas for $\tau(n)$ are given through the Mahler measure of the associated Laurent polynomial $L(z)=2k-\sum\limits_{i=1}^{k}(z^{s_k}+z^{-s_k}).$

The main results of present paper are the following. Let  $F(x)=\sum\limits_{n=1}^\infty\tau(n)x^n$ be the generating function for the number  of spanning trees $\tau(n)$  for the family of circulant graphs $C_{n}(s_1,s_2,\ldots,s_k)$ or  $C_{2n}(s_1,s_2,\ldots,s_k,n).$
We  show  (see Theorem~\ref{theoremR1} and  Theorem~\ref{theoremR2} respectively) that
$F(x)$ is a rational function with integer coefficients satisfying the property $F(x)=F(1/x).$  The obtained results are illustrated  by a series of examples.

The idea to write this paper was born during the discussion with professor Sergei Lando on the International Conference and PhD-Master Summer School on Graphs and Groups, Spectra and Symmetries (G2S2)  held on August 15 - 28, 2018, in Novosibirsk, Akademgorodok, Russia.

\section{Basic definitions and preliminary facts}

Let $s_1,s_2,\ldots,s_k$ be integers such that $1\leq s_1<s_2<\ldots<s_k\leq\frac{n}{2}.$
The graph $C_{n}(s_1,s_2,\ldots,s_k)$ with $n$ vertices $0,1,2,\ldots,~{n-1}$ is called
\textit{circulant graph} if the vertex $i,\, 0\leq i\leq n-1$ is adjacent to the vertices
$i\pm s_1,i\pm s_2,\ldots,i\pm s_k\ (\textrm{mod}\ n).$ When $s_k<\frac{n}{2}$ all vertices of
a graph have even degree $2k.$ If $n$ is even and $s_k=\frac{n}{2},$ then all vertices have odd degree $2k-1.$
It is well known  that the circulant $C_n(s_1,s_2,\ldots,s_k)$ is connected if and
only if ${\rm gcd}\,(s_1,s_2,\ldots,s_k,n)~=~1. $  More generally, the number of connected
components of $C_n(s_1,s_2,\ldots,s_k)$ is $d=\textrm{gcd}\,(s_1,s_2,\ldots,s_k,n),$
with each of the vertices $0,1,...,d - 1$ lying in different components, and with each
component being isomorphic to $C_{n/d}(s_1/d,s_2/d,\ldots,s_k/d).$ So, for $d>1$ graph
is disconnected and has no spanning trees. In what follows, all graphs are supposed to
be connected.

Two circulant graphs $C_{n}(s_1,s_2,\ldots,s_k)$ and $C_{n}(\tilde{s}_1,\tilde{s}_2,\ldots,\tilde{s}_k)$
of the same order are said to be conjugate by multiplier if there exists an integer $r$
coprime to $n$ such that $\{\tilde{s}_1,\tilde{s}_2,...,\tilde{s}_k\}=\{rs_1,rs_2,\ldots,rs_k\}$
as subsets of $\mathbb{Z}_n.$ In this case, the graphs are isomorphic, with multiplication
by the unit $ r\,(\textrm{mod}\,n)$  giving an isomorphism.

In 1967, A.~\'Ad\'am conjectured that two circulant graphs are isomorphic if and only if they are
conjugate by a multiplier \cite{Adam}.  The following example shows that the \'Ad\'am Conjecture is not true. The graphs $C_{16}(1,2,7)$ and $C_{16}(2,3,5)$ are isomorphic, but they are not conjugate by a multiplier \cite{CondGra}. A complete solution of the isomorphism problem for circulant graphs was obtained by M.~Muzychuk  \cite{Muz}.

During the paper, we will use the  basic properties of Chebyshev polynomials. Let $T_n(z)=\cos(n\arccos z)$ and $U_{n-1}(z)={\sin(n\arccos z)}/{\sin(\arccos z)}$ be the Chebyshev polynomials of the first and second kind respectively.

Then $T^{\prime}_n(z)=n\,U_{n-1}(z),\,T_n(1)=1,\,U_{n-1}(1)=n.$ For $z\neq0$ we have $T_n(\frac12(z+z^{-1}))=\frac12(z^n+z^{-n}).$

See monograph \cite{MasHand} for more advanced properties.

\section{Complexity of circulant graphs of even valency}\label{count}

The aim of this section is to find a new formula for the numbers of spanning trees of
circulant graph $C_{n}(s_1,s_2,\ldots,s_k).$ It will be based on our earlier results \cite{MedMed3,MedMed2018}, where the numbers of spanning trees was given in terms of the Chebyshev polynomials.

By Theorem~1, formula (4) from \cite{MedMed2018}, we have the following result.
\bigskip

\begin{theorem}\label{theorem1}
The number of spanning trees $\tau(n)$ in the circulant graph  $C_{n}(s_1,s_2,\ldots,s_k),$
$1\le s_1< s_2<\ldots< s_k<\frac{n}{2},$ is given by the formula
\begin{equation}\label{oldform}\tau(n)= \frac{(-1)^{n (s_k-1)} n }{q}\prod_{p=1}^{s_k-1}(2T_n(w_p)-2),\end{equation} where $q=s_1^2+s_2^2+\ldots+s_k^2,\,w_p,\,p=1,2,\ldots,s_k-1$
are all the roots of the algebraic equation $P(w)=0,$ where $$P(w)= \sum_{j=1}^k \frac{T_{s_j}(w)-1}{w-1}$$
and $T_k(w)$ is the Chebyshev polynomial of the first kind.
\end{theorem}
 We use the following elementary lemma.
 \bigskip

 \begin{lemma}\label{element}  Let  $T_n(w)$  be the Chebyshev polynomial of the first kind and $w=\frac12(z+\frac{1}{z}).$  Then
 $$T_n(w)=\frac12(z^n+\frac{1}{z^n}).$$
 \end{lemma}
 \bigskip

 By making use of substitutions $w=\frac12(z+\frac{1}{z})$  and $w_p=\frac12(z_p+\frac{1}{z_p}),$  by Lemma~\ref{element}   deduce that $2T_n(w_p)-2=-(z_p^n-1)(z_p^{-n}-1)$  and $P(w)=\sum_{j=1}^k\frac{(z^{s_j}-1)(z^{-s_j}-1)}{(z-1)(z^{-1}-1)}.$  We set $Q(z)=  \sum_{j=1}^k (z^{s_j}-1)(z^{-s_j}-1)$  and note that $Q(1)=Q^{\prime}(1)=0$  and $Q^{\prime\prime}(1)=-2(s_1^2+s_2^2+\ldots+s_k^2)=-2q<0.$  Hence, $P(w)=0$  if and only if  $w=\frac12(z+\frac{1}{z}),$  where $z$ is different from $1$  root of the equation $Q(z)=0.$

Now, Theorem~\ref{theorem1} can be restated in the following way.
\bigskip

\begin{theorem}\label{theorem2}
The number of spanning trees $\tau(n)$ in the circulant graph $C_{n}(s_1,s_2,\ldots,s_k),$
$1\le s_1< s_2<\ldots< s_k<\frac{n}{2}$ is given by the formula
\begin{equation}\label{oldform}\tau(n)= \frac{(-1)^{(n+1) (s_k-1)} n }{q}\prod_{p=1}^{s_k-1}(z_p^n-1)(z_p^{-n}-1),\end{equation} where $z_p,\,z_p^{-1},\, p=1,2,\ldots,s_k-1$
are all the roots different from $1$ of the   equation $Q(z)=0,$ and
$$Q(z)=  \sum_{j=1}^k (z^{s_j}-1)(z^{-s_j}-1).$$
\end{theorem}

    \section{Generating function for the circulant graphs    of even valency.}

The main result  of this section is the following theorem.
\bigskip
\begin{theorem}\label{theoremR1}Let  $\tau(n)$  be the number of spanning trees in the circulant graph $C_{n}(s_1,s_2,\ldots,s_k)$ of even valency. Then
$$F(x)=\sum\limits_{n=1}^\infty\tau(n)x^n$$  is  a rational  function with integer coefficients. Moreover, $F(x)=F(1/x).$
\end{theorem}

\bigskip
The proof of  the theorem is based on Theorem~\ref{theorem2} and the following proposition.
\bigskip

\begin{proposition}\label{proposition1} Let $R(z)$  be a degree $2s$ polynomial with integer coefficients.  Suppose that all  the roots  of the polynomial  $R(z)$  are $\xi_1,\xi_2,\ldots,\xi_{2s-1},\xi_{2s}.$   Then $$F(x)=\sum\limits_{n=1}^\infty n\prod\limits_{j=1}^{2s}(\xi_{j}^{n}-1)x^n$$  is  a rational  function with integer coefficients.

Moreover, if  $\xi_{j+s}=\xi_{j}^{-1},\,j=1,2,\ldots,s,$  then $F(x)=F(1/x).$
\end{proposition}

\begin{pf} First of all, we note that $F(x)=x\frac{d G(x)}{d x},$
where $$G(x)=\sum\limits_{n=1}^\infty\prod\limits_{j=1}^{2s}(\xi_{j}^{n}-1)x^n.$$

Denote by $\sigma_k=\sigma_k(x_1,x_2,\ldots,x_{2s})$ the $k$-th basic symmetric polynomial in variables $x_1,x_2,\ldots,x_{2s}.$  Namely,  $$\sigma_0=1,\,\sigma_1=x_1+x_2+\ldots+x_{2s}, \,\sigma_2=x_1x_2+x_1x_3+\ldots+x_{2s-1}x_{2s},\ldots, \sigma_{2s}=x_1x_2\ldots x_{2s}.$$  Then
$$G(x)=G_{2s}(x)-G_{2s-1}(x)+\ldots-G_1(x)+ G_0(x),$$  where $$G_k(x)= \sum\limits_{n=1}^\infty\sigma_k(\xi_1^n,\xi_2^n,\ldots, \xi_{2s}^n)x^n,\,k=0,1,\ldots,2s.$$  We have
$\sigma_k(\xi_1^n,\xi_2^n,\ldots, \xi_{2s}^n)=\sum\limits_{\substack{1\le{j_1}<{j_2}<\ldots<{j_k}\le2s}}\xi_{j_1}^n\xi_{j_2}^n\ldots\xi_{j_k}^n.$  Hence,

$$G_k(x)=\sum\limits_{\substack{1\le{j_1}<{j_2}<\ldots<{j_k}\le2s}}\frac{\xi_{j_1}\xi_{j_2}\ldots\xi_{j_k}}{1-\xi_{j_1}\xi_{j_2}\ldots\xi_{j_k} x}$$   and \begin{equation}\label{derivative}
F_k(x)=x\frac{d G_k(x)}{d x}=
\sum\limits_{\substack{1\le{j_1}<{j_2}<\ldots<{j_k}\le2s}}\frac{ \xi_{j_1}\xi_{j_2}\ldots\xi_{j_k} x}{(1-\xi_{j_1}\xi_{j_2}\ldots\xi_{j_k} x)^2}.
\end{equation}  We note that $F_k(x)$ is a   symmetric function in the roots $\xi_1,\xi_2,\ldots,\xi_{2s}$  of the integer  polynomial $R(x).$    By the Vieta theorem, $F_k(x)$ is   a rational function with integer coefficients. Since

\begin{equation}\label{alternative} F(x)=F_{2s}(x)-F_{2s-1}(x)+\ldots-F_1(x)+ F_0(x),\end{equation}
the same is true for $F(x).$

To prove the second statement of the proposition,  consider the product $\xi=\xi_{j_1}\xi_{j_2}\ldots\xi_{j_k}.$
Since $\xi_{j+s}=\xi_{j}^{-1},$   the term $\frac{\xi x}{(1-\xi x)^2}$ comes into  the sum   (\ref{derivative}) together with
$\frac{ \xi^{-1} x}{(1-\xi^{-1}x)^2}.$ One can check   that  the function  $\varphi(x)=\frac{\xi x}{(1-\xi x)^2}+\frac{\xi^{-1}x}{(1-\xi^{-1}x)^2} $ satisfy the condition $\varphi(x)=\varphi(1/x).$  Hence, for any   $k=0,1,\ldots, 2s$  we have $F_k(x)=F_k(1/x).$ By (\ref{alternative}), we finally obtain $F(x)=F(1/x).$
\end{pf}

{\bf Proof of Theorem~\ref{theoremR1}.} We employ Proposition~\ref{proposition1} to prove the theorem. To do this, we consider the polynomial $R(z)=z^{2s_k}Q(z)/(z-1)^2.$  Note that $R(z)$ is an integer polynomial of order $2s=2s_k-2.$ Recall that $Q(1)=Q^{\prime}(1)=0$  and $Q^{\prime\prime}(1)=-2q<0.$ Hence, all the roots of the polynomial $R(z)$ are the roots of  $Q(z)$ different from $1.$ Because  of an evident property $R(z)=R(\frac{1}{z}),$ the polynomial $R(z)$ satisfies the conditions of  Proposition~\ref{proposition1}. By Theorem~\ref{theorem2},  the generating function $F(x)=\sum_{n=1}^\infty\tau(n)x^n$  can be represented in the form
$$F(x)=\frac{(-1)^{s_k-1}}{q}\sum\limits_{n=1}^\infty n\prod_{p=1}^{s_k-1}(z_p^n-1)(z_p^{-n}-1)((-1)^{s_k-1}x)^n,$$
where $z_p,\,z_p^{-1},\, p=1,2,\ldots,s_k-1$
are all the roots  of the   polynomial $R(z).$  By Proposition~\ref{proposition1}, $F(x)$  is a rational function with integer coefficients satisfying $F(x)=F(\frac{1}{x}).$
\section{Complexity of circulant graphs of odd valency}\label{oddcomplexity}

The aim of this section is establish the rationality of generating function for the numbers of spanning trees in the circulant graph $C_{2n}(s_1,s_2,\ldots,s_k,n)$ of odd valency. The following result was obtained in (\cite{MedMed2018}, Theorem~2).

\bigskip

\begin{theorem}\label{odddegree2}
Let $C_{2n}(s_1,s_2,\ldots,s_k,n),\,1\leq s_1<s_2<\ldots<s_k<n,$ be a circulant graph of odd
valency. Then the number $\tau(n)$ of spanning trees in the graph $C_{2n}(s_1,s_2,\ldots,s_k,n)$
is given by the formula
$$\tau(n)= \frac{n\,4^{s_k-1}}{q}\prod_{p=1}^{s_k-1}(T_n(u_p)-1)\prod_{p=1}^{s_k}(T_n(v_p)+1),$$
where $q=s_1^2+s_2^2+\ldots+s_k^2,$ the numbers $u_p,\,p=1,2,\ldots,s_k-1\text{ and } v_p,\,p=1,2,\ldots,s_k$
are respectively the roots of the algebraic equations $P(u)-1=0,\,u\ne1$ and $P(v)+1=0,$
where $P(w)=2k+1-2\sum\limits_{i=1}^{k}T_{s_i}(w)$ and $T_k(w)$ is the Chebyshev polynomial
of the first kind.
\end{theorem}

Taking into account the identity $T_n(z +\frac{1}{z})=\frac12(z^n+\frac{1}{z^n}),$ we rewrite Theorem~\ref{odddegree2} in the following way.
\bigskip

\begin{theorem}\label{theorem22}
The number of spanning trees $\tau(n)$ in the circulant graph $C_{2n}(s_1,s_2,\ldots,s_k,n),\,1\leq s_1<s_2<\ldots<s_k<n$ of odd valency, is given by the formula
\begin{equation}\label{oldformore}
\tau(n)= \frac{(-1)^{s_{k}-1}n}{2q}\prod_{p=1}^{s_k-1}(z_p^n-1)(z_p^{-n}-1)\prod_{p=1}^{s_k}(\zeta_p^n+1)(\zeta_p^{-n}+1),
\end{equation}
where $z_p,\,z_p^{-1},\, p=1,2,\ldots,s_k-1$
are all the roots different from $1$ of the   equation $Q(z)=0,\, \zeta_p,\,\zeta_p^{-1},\, p=1,2,\ldots,s_k$
are all the roots of the   equation $Q(z)+2=0,$ and
$$Q(z)=  \sum_{j=1}^k (z^{s_j}-1)(z^{-s_j}-1).$$

\end{theorem}

Using this result and repeating the arguments from the proofs of Proposition~\ref{proposition1} and Theorem~\ref{theoremR1}, we obtain the following theorem.
\bigskip

\begin{theorem}\label{theoremR2}Let  $\tau(n)$  be the number of spanning trees in the circulant graph $C_{2n}(s_1,s_2,\ldots,s_k,n)$ of odd valency. Then
$$F(x)=\sum\limits_{n=1}^\infty\tau(n)x^n$$  is  a rational  function with integer coefficients. Moreover, $F(x)=F(1/x).$
\end{theorem}

\section{Examples}\label{tables}

Everywhere in the section below we use the following notation $w=\frac{1}{2}(x+\frac{1}{x}).$  Since the generating function
$F(x)=\sum\limits_{n=1}^\infty\tau(n)x^n$ satisfies $F(x)=F(1/x),$  it is also a rational function
 in  $w.$
\begin{enumerate}
\item{\textbf{Graph $C_n(1,2)$}.}
By \cite{BoePro} and \cite{MedMed2018}   we have $\tau(n)=nF_n^2,$ where $F_n$ is the $n$-th Fibonacci number.  Hence,
$$\sum\limits_{n=1}^\infty\tau(n)x^n=\frac{1 - 2 w + 2 w^2}{(1 + w) (-3 + 2 w)^2}.$$

\item{\textbf{Graph $C_n(1,3)$}.}
By \cite{MedMed2018}, Example~2, we have
$$\tau(n)=\frac{2n}{5}(T_n(-\frac12-\frac{i}{2})-1)(T_n(-\frac12+\frac{i}{2})-1).$$
Hence,
$$\sum\limits_{n=1}^\infty\tau(n)x^n=\frac{(1+w)(1-w-2w^2+11w^3+8w^4-16w^5+4w^7)}
{2(-1+w)(-1-3w-3w^2+2w^4)^2}.$$

\item{\textbf{Graph  $C_n(2,3)$}.}
In this case, by \cite{MedMed2018}, Example~3,
$\tau(n)=\frac{4n}{13}(T_n(\theta_1)-1)(T_n(\theta_2)-1),$ where $\theta_{1,2}=\frac{1}{4}(-3\pm i\sqrt{3}).$ Now
$$\sum\limits_{n=1}^\infty\tau(n)x^n=\frac{2(15+27w-17w^2-28w^3+216w^4+64w^5-224w^6+64 w^8)}
{(-1+w)(-3-24w-28w^2+16w^4)^2}.$$

\item{\textbf{Graph M\"obius ladder $C_{2n}(1,n).$}}
In this case, by \cite{BoePro} and Theorem~\ref{odddegree2} we
have $\tau(n)=n(T_n(2)+1)$ and
$$\sum\limits_{n=1}^\infty\tau(n)x^n=\frac{5 - 7 w + 3 w^2}{2 (-2 + w)^2 (-1 + w)}.$$

\item{\textbf{Graph $C_{2n}(1,2,n)$}.} By Theorem~\ref{odddegree2} we obtain
$$\tau(n)=\frac{4 n}{5} (-1 + T_n(-3/2)) (1 + T_n (\theta_1)) (1 + T_n (\theta_2)),$$
where $\theta_{1,2}=\frac14(-1\pm\sqrt{33}).$  As a result,
$$\sum\limits_{n=1}^\infty\tau(n)x^n=\frac{G(w)}{H(w)},$$ where

\begin{eqnarray*}
G(w)&=&- 422137930752 + 83760842880 w + 3098738380752 w^2 \\
&&- 883496472336 w^3 - 9652661216757 w^4 +3106365668883 w^5 \\
&&+ 17117057275452 w^6 - 5625024610977 w^7 - 19349031968433 w^8 \\
&&+ 6096193861050 w^9 + 14674263153534 w^{10} - 4191031652328 w^{11} \\
&&- 7620826953696 w^{12} + 1847068052640 w^{13} + 2715353668608 w^{14} \\
&&- 509714698368 w^{15} - 654968016896 w^{16} + 82874686976 w^{17} \\
&&+ 103379897344 w^{18} - 6793590784 w^{19} - 9988546560 w^{20} \\
&&+ 70311936 w^{21} + 531562496 w^{22} + 34570240 w^{23} \\
&&- 14614528 w^{24} - 2097152 w^{25} + 393216 w^{26}
\end{eqnarray*}

and
\begin{eqnarray*}
H(w)&=&(-1 + w) \big((3 + 2 w) (-4 + w + 2 w^2) (13 + 16 w + 4 w^2)  \\
&&(96 - 216 w + 165 w^2 - 48 w^3 + 4 w^4) (87 + 9 w - 90 w^2 - 12 w^3 + 8 w^4)\big)^2.
\end{eqnarray*}

\end{enumerate}

%%%%%%%%%%%%%%%%%%%%%%%%%%%%%%%%

 The results of this work were partially supported by the Russian Foundation for Basic Research,
(grants 18-01-00420 and 18-501-51021).

\end{document}